\documentclass[a4paper,british,english]{amsart}

\usepackage{amsmath,amssymb,amsthm,amscd,epsfig}
\def\R{{\Bbb R}}

\def\N{{\Bbb N}}

\def\d{\mbox{d}}

\newcommand{\vp}{\varepsilon}

\newcommand{\bb}[1]{\mathbb{#1}}

\newcommand{\cl}[1]{\mathcal{#1}}

\theoremstyle{plain}

\newtheorem{thm}{Theorem}

\newtheorem{cor}{Corollary}

\newtheorem{pro}{Proposition}

\theoremstyle{definition}

\theoremstyle{remark}

\begin{document}

\title{Upper estimates for stable dimensions of fractal sets with variable numbers of foldings}

\author{Eugen Mihailescu and Bernd Stratmann}
\date{}

\maketitle

\begin{abstract}

For a hyperbolic map $f$ on a saddle type fractal $\Lambda$ with self-intersections, the number of $f$-preimages of a point $x$  in $\Lambda$ may depend on $x$.  This makes estimates of the stable dimensions more difficult than for diffeomorphisms or for maps which are constant-to-one.  We employ  the thermodynamic formalism in order to derive estimates for the stable Hausdorff dimension function $\delta^s$ on $\Lambda$, in the case when $f$ is conformal on local stable manifolds. These estimates are in terms of a continuous function  on $\Lambda$ which bounds the preimage counting function from below. As a corollary we obtain that if  $\delta^s$ attains its maximal possible value in $\Lambda$, then the stable dimension is constant throughout $\Lambda$, whereas the preimage counting function is constant on at least an open and dense subset of $\Lambda$. In particular, this shows that  if at some point in $\Lambda$ the stable dimension is equal to the analogue of the similarity dimension in the stable direction at that point, then $f$
  behaves very much like a homeomorphism on $\Lambda$. Finally, we also obtain results about the stable upper box dimension for these type of fractals. We end the paper with a discussion of two explicit examples.
\end{abstract}

\

\textbf{Mathematics Subject Classification 2000:} 37D35, 37F15
37D45, 37F10.

\

\textbf{Keywords:}   Non-invertible hyperbolic maps; thermodynamic formalism; basic sets of saddle type; stable dimensions; conformal maps;  fractals with overlap; Hausdorff dimension; box-counting dimension.

\

\section{Introduction and statement of results.}
In this paper we investigate fractal sets $\Lambda$ of saddle type which are invariant under a non-invertible $\mathcal{C}^2$-endomorphism $f$ of a Riemann manifold $M$ into itself. These fractals are {\em basic sets of $f$}, meaning that $\Lambda$ is  compact
and  $f$-invariant such that $f|_\Lambda$ is topologically transitive and such that there exists a neighbourhood $U$ of $\Lambda$ satisfying $\Lambda = \bigcap_{n \in \mathbb Z} f^n(U)$. 
The fact that $f$ is non-invertible produces complicated overlaps and foldings within $\Lambda$, which influence the Hausdorff dimension of the sections through $\Lambda$ and the number of overlaps does not necessarily has to be constant.
We will always assume that $f$ is hyperbolic on $\Lambda$ in the sense of Ruelle \cite{Ru-carte89}, that is,  for each backward orbit  $\hat x = (x, x_{-1}, x_{-2}, \ldots)$ of $x$ in $\Lambda$, where $f(x_{-1})=x$ and $ f(x_{-(i+1)}) = x_{-i} \in \Lambda$ for all $i\in \N$,   there exists a continuous splitting of the tangent bundle over the space $\hat \Lambda$ of all backward orbits of elements of $\Lambda$, called the {\em natural extension} (or inverse limit) of the tuple $(\Lambda, f)$,  into stable spaces $E^s_x$ and unstable spaces $E^u_{\hat x}$. It is well-known that $\hat \Lambda$ is a compact metric space and that the lift $\hat f: \hat \Lambda \to \hat \Lambda$ of $f$ to $\hat \Lambda$, given by $\hat f(\hat x) := (f(x), x, x_{-1}, x_{-2}, \ldots)$, is a homeomorphism. Note that natural extensions play an important role in the study of the dynamics of endomorphisms (see for instance \cite{Ru-carte89,M-DCDS06}). As in the diffeomorphism case, for a hyperbolic endomorphism $f$ on $\Lambda$ there exist local stable manifolds $W^s_r(x)$ and local unstable manifolds $W^u_r(\hat x)$, for each $x \in  \Lambda$ and $\hat x \in \hat \Lambda$. Note that there may be  infinitely many  local unstable manifolds through a given point in $\Lambda$ and, unlike in the diffeomorphism case, these do not necessarily give rise to a foliation. 
\\
In this paper we will consider  the {\em stable dimension at $x \in \Lambda$},  which is given by
 $$\delta^s(x):= \dim_{H}(W^s_r(x) \cap \Lambda), $$ where $\dim_{H}$ refers to the Hausdorff dimension. To give estimates for the stable dimension is by far more delicate than for the unstable dimension $\delta^u(\hat x):= \dim_{H}(W^u_r(\hat x) \cap \Lambda)$.  In fact, in  \cite{M-DCDS06} it was shown that $\delta^u(\hat x)$ is constant on $\hat \Lambda$ and that its value is given by the unique zero of the pressure function $P_{\hat f|_{\hat \Lambda}}(-t\log |Df_u|)$, where $|Df_u(x)|$
denotes the norm of the derivative of $f$ restricted to $E^u_{\hat x}$. However, for the stable dimension we can not expect that a similar formula holds in general. 
\\ Before we state our main result, let us point out that in this paper we consider a special type of hyperbolic endomorphisms which will be called \textit{c-hyperbolic}. A map $f$  is {\em c-hyperbolic} on $\Lambda$ if it is hyperbolic as an endomorphism over $\Lambda$,  if it is conformal on all local stable manifolds and if $\Lambda$ does not contain any critical points of $f$.\\
Also, let us introduce the {\em preimage counting function} $\Delta:\Lambda \to \N,$ which is given for each $x\in \Lambda$ by $$\Delta(x) := \hbox{Card} \left(f^{-1}(x) \cap \Lambda\right)$$
One immediately verifies that  $\Delta$  is upper semi-continuous and bounded  on $\Lambda$ (see e.g. \cite[Lemma 1]{MU-BLMS}).
Moreover,  the {\em stable potential function} $\Phi^s$ on $\Lambda$ is defined by  $\Phi^s(x) := \log |Df_s(x)|$, where $|Df_s(x)|$
denotes the norm of the derivative of $f$ restricted to $E^s_x$. 
We are now in the position to state the main result of this paper.
\begin{thm} \label{upper}
Let $f:M \to M$ be a $\mathcal{C}^2$-endomorphism which is c-hyperbolic on a basic set $\Lambda$ of $f$ and for which there exists a continuous function $\omega:\Lambda \to \R$ such that $\Delta(x) \ge \omega(x)$, for all $ x \in \Lambda$. It then follows that
\[\delta^s(x) \le t_\omega, \]
where $t_\omega$ refers to the unique zero of the pressure function $t \mapsto P(t \Phi^s- \log \omega)$ 
 associated with the potential function $t \Phi^s- \log \omega$.
\end{thm}


Let us point out that one of the difficulties in proving this theorem is that the map $f$ is not necessarily expanding and that its inverse branches do not necessarily contract small balls. In fact, some directions may be even expanding in backward time. Another difficulty is that the number of preimages of a point that remain in $\Lambda$ is not always constant. 


The reader might like to recall that in their pioneering work Bowen \cite{Bo-dim} and Ruelle \cite{Ru-dim} employed  the thermodynamic formalism in order to derive dimension formulae for rational maps. In fact, in the diffeomorphism case, it turned out that the stable and the unstable dimension can in general be  computed both as the zero of the pressure function of the stable potential, respectively the unstable potential (see \cite{MM}); (for further applications of the thermodynamic formalism in dimension theory, we refer to \cite{Ba}, \cite{P}).
It is important to note that  for an endomorphism $f$ in higher dimension, a hyperbolic basic set  is not necessarily totally invariant. This is of course significantly different from the case of Julia sets of rational maps in the complex one dimensional case.
Examples of perturbations of toral endomorphisms which are Anosov and whose unstable manifolds depend on the whole prehistory were given in \cite{Pr}. 
Another class of non-invertible hyperbolic maps with crossed invariant horseshoes was given by Bothe in \cite{Bot}. 
Also, Simon \cite{Si} gave another class of non-invertible endomorphisms,  for which the Hausdorff dimension of the associated attractors can be computed with the help of a pressure formula just as in the invertible case.  Examples of non-linear hyperbolic skew products of Cantor sets with overlaps in their fibres were given in \cite{M-MZ}, where 
the strongly non-invertible character of these maps has been established, and where it was shown that these skew products are far away from being constant-to-one. In \cite{M-MZ} it was also shown that there are maps  for which there exist Cantor sets in each of their fibres, such that through each point of these sets there pass uncountably many different local unstable manifolds.  
Also let us mention that  yet another class of c-hyperbolic endomorphisms can be found by considering hyperbolic basic sets  of saddle type, for holomorphic maps $f: \bb P^2\bb C \to \bb P^2\bb C$ on the 2-dimensional complex projective space (\cite{M-DCDS06}).

The paper continues by showing that  an application of Theorem \ref{upper} gives rise to the following proposition. In here, we consider the situation in which $\delta^{s}$  attains a maximal value and show that in this case, $\delta^{s}$ has to be constant throughout $\Lambda$ and that $\Delta$ has to be equal to its least value $d$ on an open dense subset.
\begin{pro}\label{max}
If in addition to the assumptions in Theorem \ref{upper} we have that the minimal value of $\Delta$ on $\Lambda$ is equal to $d$ and that there exists a point $x \in \Lambda$ at which $\delta^s$ is equal to the unique zero $t_d$ of the pressure function $t \mapsto P(t\Phi^s - \log d)$, then $\Delta$ is equal to $d$ on an open dense subset of $\Lambda$ and  $\delta^s(y)$ is equal to $  t_d$, for all $y \in \Lambda$.
\end{pro}
Note that the latter proposition can be  applied  in particular in the case in which $d$ is equal to $1$ and where there is no overlap. In this situation the stable dimension is equal to the similarity dimension, and the proposition then guarantees  that there exists an open dense set 
of points in $\Lambda$ at which $f$ has precisely one preimage in $\Lambda$. Therefore, in this case  the map behaves almost like a homeomorphism when restricted to $\Lambda$. This particular situation is somewhat  parallel to a result of Schief \cite{S}, although the setting and proofs are completely different. We summarise these results in the following corollary. 
\begin{cor}\label{inj}
Let $f:M \to M$ be a $\mathcal{C}^2$-endomorphism which is c-hyperbolic on a basic set $\Lambda$ of $f$ and for which there
 exists a point $x \in \Lambda$ such that $\delta^s(x)$ is equal to the unique zero $t_1$ of the pressure function $t \mapsto P(t\Phi^s)$. Then there exists an open dense set of points in $\Lambda$ at which $f$ has precisely one preimage in $\Lambda$. Moreover,  we have that  $\delta^s(y) = t_1$, for all $y \in \Lambda$.
\end{cor}  

Also, in \textbf{Corollary \ref{horseover}} from Section 4 we will show how our  Corollary \ref{inj}  can be applied to a class of translations of horseshoes with overlaps, previously studied by Simon and Solomyak in \cite{SS}. 

Let us now remark that a combination  of Theorem \ref{upper} with the main theorem in \cite{MU-BLMS} gives rise to the following result.
\begin{cor}\label{locconst}
If in addition to the assumptions in Theorem \ref{upper} we have that the preimage counting function $\Delta$ is locally constant on $\Lambda$, then it follows that $\delta^s(x) = t_\omega$, for all $x \in \Lambda$. Here, $t_{\omega}$ is given as in Theorem \ref{upper}.
\end{cor}
Finally, we consider the {\em stable upper box dimension} $\beta^{s}(x)$ which is given by  the upper box-counting dimension $\overline{\dim}_{B}(W^s_r(x) \cap \Lambda)$ of the set $W^s_r(x) \cap \Lambda$, for each $x \in \Lambda$. (For a general discussion of the upper box dimension for fractal sets we refer to   \cite{Mat} and \cite{P}).  We show that  this dimension function is constant  throughout $\Lambda$ and that in the situation in which $\Delta$ is bounded from below, similarly as in Theorem \ref{upper},  one derives an upper bound for its value. These results are summarised in the following proposition.
\begin{pro}\label{upperboxdim} Let $f:M \to M$ be a $\mathcal{C}^2$-endomorphism which is c-hyperbolic on a basic set $\Lambda$ of $f$. Then the following hold. \begin{itemize}
\item[(a)] If there exists a continuous function $\omega:\Lambda \to \R$ such that $\Delta(x) \ge \omega(x)$, for all $ x \in \Lambda$, then we have, with $t_{\omega}$ given as in Theorem \ref{upper},  \[ \beta^{s}(y)  \le t_\omega, \mbox{ for all } y \in \Lambda.\]
\item[(b)]   The function $\beta^{s}$  is constant on $\Lambda$. 
\end{itemize}
\end{pro}

In particular the above results apply for hyperbolic basic sets of saddle type for holomorphic maps $f: \mathbb P^2\bb C \to \bb P^2 \bb C$. 
We will end the paper by giving two further explicit examples in which  the above results can be applied. Our first example 
will be concerned with certain horseshoes with overlaps in $\mathbb R^3$ considered in \cite{SS}. The second example will be on basic sets for a family of hyperbolic skew products studied in  \cite{M-MZ}.

We close this introduction with some comments on how the results in this paper relate to previous work in this area.

In \cite{Fa} (see also \cite{So} and \cite{PS}) Falconer  studied self-affine fractals with overlaps obtained from  finitely many linear contractions $T_i(x) = \lambda_i x, i = 1,..., \ell$ in $\R$ satisfying $0 < |\lambda_i| < 1$ and $\sum_{1 \le i \le \ell} |\lambda_i| < 1$.  He showed that  
the Hausdorff dimension of the invariant set of the family of translated contractions  $\{T_i + a_i,: 1 \le i \le \ell\}$ is equal to $s$, for Lebesgue almost all  $(a_1, \ldots, a_\ell) \in   \R \times ... \times \R$; where $s$ represents the \textit{similarity dimension}, defined as the solution of the equation $$\sum_{1 \le i \le \ell} |\lambda_i|^s = 1$$
We remark that this result may be extended also to similarities on $\mathbb R^n$. 
However, the result fails if the condition 
$\sum_{1 \le i \le \ell} |\lambda_i|<1$  is not satisfied, as observed by Edgar (\cite{E}), who  based his argument on a result by Przytycki and Urba\'nski (\cite{PU-89}). Indeed, if $T_1 = T_2 =\left(\begin{array}{ll}
                                        1/2 & 0 \\
                                       0 & \lambda
  \end{array} \right)$ and if $|\lambda| > \frac 12$, then for Lebesgue almost every $a = (a_1, a_2) \in \mathbb R^2$  the attractor $\Lambda(a)$ of the system $\{T_1 + a_1, T_2 + a_2\}$ stays to be the same; and moreover if $1/\lambda$ is a Pisot  number (that is, an algebraic integer such that the absolute value of all its algebraic conjugates is less than 1), then $\dim_{H}(\Lambda(a)) < 2 - (\log(1/\lambda)) / \log 2$ (see e.g. \cite{So}).  This shows that  fractals originating from overlapping constructions can have Hausdorff dimension \textit{less} than their similarity dimension.

In \cite{S} Schief  considered self-similar fractal sets $K$ and showed that if for the similarity dimension $\sigma$ of $K$ one has that the $\sigma$-dimensional Hausdorff measure   $\mathcal{H}^\sigma(K)$ is positive,  then $K$ satisfies the strong open set condition, that is, the system behaves similar to a homeomorphism on $K$. Note that this result is in the spirit of our results in this paper, although the setting and the ideas of our proofs  differ significantly from the approach in \cite{S}. More precisely, the assumptions in Proposition \ref{max} are much weaker  than the ones in \cite{S}. Namely, in order to obtain the "almost injectivity" of the system associated with $\Lambda$, we only require that the stable dimension $\delta^s(x)$ is equal to the zero $t_1$ of the pressure function $t \to P(t \Phi^s)$, for some $x \in \Lambda$; we do not require that $\mathcal{H}^{t_1}(W^s_r(x) \cap \Lambda) >0$. 
 In our case $t_1$ is the analogue of the similarity dimension in the stable direction, in the sense that it represents the dimension which one would obtain if the system would be invertible. In particular if there exists some $x \in \Lambda$ for which $\mathcal{H}^{t_1}(W^s_r(x) \cap \Lambda) >0$ is positive,  then we have that the stable dimension is everywhere equal to $t_1$ and that there exists an open dense set of points in $\Lambda$ which have precisely one preimage in $\Lambda$.

%

Finally, in  \cite{MU-BLMS} Mihailescu and Urbanski studied  c-hyperbolic maps on $\Lambda$ for which $\Delta$ is bounded from {\em above} by a continuous map $\eta$ on $\Lambda$. In that paper the authors obtain the result that  $\delta^s(x) \ge t_\eta$ for all $x \in \Lambda$, where $t_\eta$ refers to the unique zero of the pressure function $t \to P(t\Phi^s - \log \eta)$. Note that the proof for the upper estimate in this paper is very different from the proof for the lower estimates in \cite{MU-BLMS}. However,  we can combine these two estimates, as done in Corollary \ref{locconst}, to obtain  that if the preimage counting function $\Delta$ is locally constant on $\Lambda$, then the stable dimension is equal to $t_{\Delta}$ throughout $\Lambda$. 

\

\section{Proof of Theorem \ref{upper}.}

For ease of exposition, let us first consider the situation in which $\omega$ is locally constant and takes on only two different positive integer values on $\Lambda$, namely $d_1$ on the set $V_1$ and $d_2$ on the set $V_2$. We then have that  $V_1 \cup V_2 = \Lambda$  and that $V_1$ and $V_2$ are  two disjoint compact subsets of $\Lambda$.   
Hence, there exists $\vp_0>0$ such that the distance $d(V_1, V_2)$ between $V_{1}$ and $V_{2}$ is greater than $\vp_0$.
For $x \in \Lambda$ and $n\in \N$, let $B_n(x, \vp) :=\{y \in \Lambda: d(f^i(y), f^i(x)) < \vp, 0 \le i \le n-1\}$ refer to the $n$-Bowen ball centred at $x$ of radius $\vp>0$.
Note that for $0<\vp<\vp_0$ we have that if $y \in B_n(x, \vp)$ then $f^i(y)$ and $f^i(x)$ both belong to  either  $V_{1}$ or  $V_{2}$, for each $0 \le i \le n-1$.    
Recall that $\Phi^s(x) := \log |Df_s(x)|, x \in \Lambda$.
Now, let   $t > t_\omega$ be fixed. By definition of $t_\omega$, we have that there exists $\beta>0$ such that $$P(t\Phi^s - \log \omega) < -\beta$$ Hence, by choosing  $\vp>0$ sufficiently small, there exists a constant $C>0$  such that for each  $n\in \N$ large enough, there exists a minimal $(n, \vp)$-spanning set $E_n$ for $\Lambda$ such that  
\begin{equation}\label{P}
 \mathop{\sum}\limits_{z \in E_n} (\text{diam} \, \, U_n(z))^t \cdot \frac{1}{\Delta(f(z) \cdot \ldots \cdot \Delta(f^n(z))} < C \, e^{- \beta n} < 1,
\end{equation}
where we have set  $U_n(z):= f^n(B_n(z, \vp)) \cap W^s_r(x) \cap \Lambda$. Note that in here we have used the fact that the set $U_n(z)$ is the intersection of an unstable tubular neighbourhood with the fixed stable manifold $W^s_r(x)$. Also, we used that  $|Df_s^n(z)|$ is  uniformly comparable to $ \text{diam} \, \,  U_n(z)$, which follows from the fact that $f$ is conformal on local stable manifolds.

In the sequel let us put  $W:=W^s_r(x) \cap \Lambda$. Hence, the aim is to show that $\dim_{H}(W) \le t$, for each $t > t_\omega$.  
The main idea of the proof is to extract suitable covers of $W$ out of the large set of covers which are given by taking $n$-preimages, such that at each step a different sum will be minimised. Note that we say that a point $y$ is a {\em $k$-preimage of $x$} if $f^k(y) = x$.  Each such $n$-preimage will be included in a Bowen ball of type $B_n(z, \vp)$, for some $z \in E_n$. This procedure is delicate, since at each step  the number of preimages of points belonging to $\Lambda$ varies. The idea is to consider the $k$ iterates of $n$-preimages, then to subdivide $\Lambda$ into various different parts and finally, to find suitable covers of these parts which minimise certain sums at the $k$-th level.

First, note that since $\Lambda$ is covered by the set of Bowen balls $\{B_n(z, \vp): z \in E_n\}$, it follows that  $\{U_n(z):  z \in E_n\}$ covers $W$. However,  this cover is far too rich and we will have to extract a suitable subcover. Indeed, by using a well known theorem by Besicovitch  (see for e.g. \cite{Mat}), there exists a subcover $\{5 U_n(z): z \in \mathcal{G}(0)\}$ of $W$ such that 
$\{U_n(z): z \in \mathcal{G}(0)\}$ consists of pairwise disjoint sets. (Note that, since $f$ is conformal on local stable manifolds, we can assume that the sets $U_n(z)$ are in fact balls, and we denote the radii of these by $r(n, z)$;  also, we write $5U_n(z)$  to denote the ball of radius $5 r(n, z)$ centred at the centre of $U_n(z)$). The next step is  to "inflate" this cover, that is,  to enlarge it to a ''richer'' cover of $W$. For this, we consider an $(n-1)$-preimage of $w$ in $\Lambda$ which we denote by $w(n-1)$, for each point $w \in W$. Let us assume that $w(n-1)\in V_1$ and hence, that  $w(n-1)$ has at least $d_1$ 1-preimages in $\Lambda$.
 Now, since $E_n$ is $(n, \vp)$-spanning, for each point $\xi \in \Lambda$, there exists at least one point $y \in E_n$ such that $\xi \in B_n(y, \vp)$. However, we cannot have two 1-preimages of some $w(n-1)$ belonging to different Bowen balls $B_n(y, \vp)$ and $B_n(y', \vp)$  such that $y$ and $y'$ are both in $\mathcal{G}(0)$. This is an immediate consequence of the fact that  $\{U_n(z): z \in \mathcal{G}(0)\}$ consists of pairwise disjoint sets.
 
  Therefore, by way of successive eliminations, we can find  $d_1$ pairwise disjoint families, denoted by  $\cl{F}(1, d_1; 1), \ldots, \cl{F}(1, d_1; d_1)$, such that $\{5 U_n(z): z \in \cl{F}(1, d_2; i)\}$ is a cover of the set $\{w \in W: w(n-1) \in V_1\}$, for each $1\leq i \leq d_{1}$. Obviously, for $w(n-1)\in V_2$ we can proceed in a similar way, which then gives rise 
to $d_2$ mutually disjoint families $\cl{F}(1, d_2; 1), \ldots, \cl{F}(1, d_2; d_2)$ for which we have that  $\{5 U_n(z): z \in \cl{F}(1, d_2; j)\}$  is a cover of $\{w \in W: w(n-1) \in V_2\}$, for each $1 \le j \le d_2$. Note that, since $d(V_1, V_2) >0$, we have that $\cl{F}(1, d_1; i) \cap \cl{F}(1, d_2; j) = \emptyset$, for all $ i$ and $j$, and that by construction we have  that the so obtained disjoint families are all contained in $E_n$.
Next, define
 $$\cl{F}(1):= \mathop{\cup}\limits_{i = 1}\limits^2\mathop{\cup}\limits_{1 \le j \le d_i}\cl{F}(1, d_i, j)$$
 and let $\cl{G}(1, d_k)$ be given, for $k\in \{1,2\}$,  by
$$ \mathop{\sum}\limits_{z \in \cl{G}(1, d_k)} \frac{(\text{diam} \,U_n(z))^t}{\Delta(f^2 (z)) \ldots \Delta(f^n (z))}= \min \left\{\mathop{\sum}\limits_{z \in \cl{F}(1, d_k; i)} \frac{(\text{diam} \,U_n(z))^t}{\Delta(f^2 (z)) \ldots \Delta(f^n (z))}:i \in\{1,...,d_{k}\} \right\}.$$
For $\cl{G}(1):= \cl{G}(1, d_1) \cup \cl{G}(1, d_2)$, we then obtain, by adding the sums over $\cl{G}(1, d_1)$ and $\cl{G}(1, d_2)$, 
\begin{equation}\label{step1}
\mathop{\sum}\limits_{z \in \cl{G}(1)} \frac{(\text{diam} \, U_n(z))^t}{\Delta(f^2(z)) \ldots \Delta(f^n (z))} \le \mathop{\sum}\limits_{z \in \mathcal{F}(1)}\frac{(\text{diam} \, U_n(z))^t}{\Delta(f(z)) \ldots \Delta(f^n(z))}.
\end{equation}
Note that here we have used the trivial fact that for each $x \in \Lambda$ we have that
$ \sum_{y \in \Lambda, f(y)=x} 1/\Delta(x) =1$.
Also, note that  the sum  over the family $\cl{G}(1)$ on the left hand side of the inequality in (\ref{step1}) is smaller than the sum over the  larger family $\cl{F}(1)$ on the right hand side. However, and this is the crucial point,  the summands on the right hand side have one more factor in their denominator than the summands on the left hand side. 

Let us now bring the argument to its next level by  enlarging the family $\cl{F}(1)$ as follows. Recall that for each $w \in W$ we have fixed an $(n-1)$-preimage $w(n-1) \in \Lambda$. We now define $w(n-2):= f(w(n-1))$ and consider not only $w(n-1)$ but  also the other 1-preimages of $w(n-2)$ in $\Lambda$. Subsequently, we will then take the 1-preimages of these 1-preimages of $w(n-2)$ and obtain new covers of $W$. 
Indeed similarly as before, if $w(n-2) \in V_1$ then we can construct, by succesive eliminations, pairwise disjoint families $\cl{F}(2, d_1; 1), ..., \cl{F}(2, d_1; d_{1}) $ by selecting the 1-preimages of the $i$-th preimage of $w(n-2)$, for each $1 \le i \le d_1$. In fact one of these families is $\cl{F}(1)$. As in the first step, the sets $\{5 U_n(z): z \in \cl{F}(2, d_1; i)\}$  cover $\{w\in W: w(n-2) \in V_1\}$, for each $i$.
Let us remark that the procedure of successive elimination works, since if we take for instance the family $\cl{F}(2, d_1; 1)$, then for an arbitrary $w \in W$ we cannot have two 1-preimages $y$ and $y'$ of $w(n-2)$ and 1-preimages $\xi$ of $y$ and $\xi'$ of $y'$ such that $\xi$ and $\xi'$ are both contained  in either $B_n(z, \vp)$ or $ B_n(z', \vp)$, for some $z, z' \in \cl{F}(2, d_1; 1)$. Indeed, since $f^2(B_n(z, \vp)) \cap f^2(B_n(z', \vp)) \ne \emptyset$, in this situation it would follow that $U_n(z) \cap U_n(z') \ne \emptyset$ and hence we would have a contradiction. This implies that there exist $d_1$ disjoint families $\cl{F}(2, d_1; i)$ corresponding to the $d_1$ 1-preimages of $w(n-2) \in V_1$. 

Clearly, we can proceed analogously in the case in which $w(n-2) \in V_2$, which then gives rise to pairwise disjoint families $\cl{F}(2, d_2; 1), ..., \cl{F}(2, d_2; d_{2}) $ for which $\{5 U_n(z): z \in \cl{F}(2, d_2; j)\}$  covers $\{w\in W: w(n-2) \in V_2\}$, for each $j$.  Note that we cannot have repetitions of points from $E_n$ when taking the union of the collections $\cl{F}(2, d_i; j)$ over all $ i \in \{1,2\}$ and $ 1 \le j \le d_i$. Indeed, if we would have two 1-preimages $y, y' \in \Lambda$ of some $w(n-2)$ and two 1-preimages $\xi, \xi' \in \Lambda$ of $y$, and $y'$ respectively, so that $\xi \in B_n(z, \vp)$ and $\xi' \in B_n(z', \vp)$, for some $z, z' \in \mathcal{F}(2, d_1; i)$, then it would follow that $U_n(z)\cap U_n(z') \ne \emptyset$,  which gives a contradiction. Moreover, by construction we have that  $\cl{F}(2, d_1; i) \cap \cl{F}(2, d_2; j) = \emptyset$, for all $ i$ and $j$. 
This follows,  since if $f^2(z) \in V_1$, for some $z \in \cl{F}(2, d_1; i)$, and if at the same time  $f^2(z') \in V_2$, for some $z' \in \cl{F}(2, d_2; j)$, then it would follow that  $V_1 \cap V_2 \neq \emptyset$ and hence, we would get a contradiction.

Now, as in the first step, for each $i \in \{1, 2\}$ and  $1 \le j \le d_i$ there exists a family  $\cl{G}(2, d_i,j)$ in  $\{\cl{F}(2, d_k; \ell):
k \in \{1,2\},1 \leq  \ell \leq d_{j}\}$ satisfying
\begin{equation}\label{step2}
\mathop{\sum}\limits_{z \in \cl{G}(2, d_i,j)} \frac{(\text{diam} \, U_n(z))^t}{\Delta(f^2(z)) \ldots \Delta(f^n(z))} \le \mathop{\sum}\limits_{z \in \cl{F}(2, d_i; j)} \frac{(\text{diam} \, U_n(z))^t}{\Delta(f(z)) \Delta(f^2(z)) \ldots \Delta(f^n(z))}.
\end{equation}
Among these so obtained families $\cl{G}(2, d_i; j)$ we now choose for each $i \in \{1,2\}$ a particular family, which will be denoted by $\cl{G}(2, d_i)$, such that  
\begin{equation}\label{step2,1}\mathop{\sum}\limits_{z \in \cl{G}(2, d_i)} \frac{(\text{diam} \, U_n(z))^t}{\Delta(f^3 (z)) \ldots \Delta(f^n (z))}=
 \min \left\{\mathop{\sum}\limits_{z \in \cl{G}(2, d_i;j)} \frac{(\text{diam} \, U_n(z))^t}{\Delta(f^3 (z)) \ldots \Delta(f^n (z))}:j \in \{1,...,d_{i}\} \right\}.\end{equation}
Combining (\ref{step2}) and (\ref{step2,1}), we now have for each $i \in \{1,2\}$ that
$$
\aligned
\sum_{z \in \cl{G}(2, d_i)} \frac{(\text{diam} \, U_n(z))^t}{\Delta(f^3 (z)) \ldots \Delta(f^n(z))} &\le \sum_{z \in \bigcup_{1 \le j \le d_i} \cl{G}(2, d_i; j)} \frac{(\text{diam} \, U_n(z))^t}{\Delta(f^2 (z)) \ldots \Delta(f^n(z))} \\
&\le \sum_{z \in \bigcup_{1 \le j \le d_i} \cl{F}(2, d_i; j)} \frac{(\text{diam} \, U_n(z))^t}{\Delta(f(z)) \Delta(f^2 (z)) \ldots \Delta(f^n (z))}.
\endaligned
$$
Therefore, by defining $$\cl{F}(2):= \bigcup_{i\in \{ 1, 2\}} \bigcup_{1 \le j \le d_i} \cl{F}(2, d_i; j) \ \text{  and  } \ \cl{G}(2):= \cl{G}(2, d_1) \cup \cl{G}(2, d_2),$$ 
we have now shown that
\begin{equation}\label{s2}
\sum_{z \in \cl{G}(2)} \frac{(\text{diam} \, U_n(z))^t}{\Delta(f^3(z)) \ldots \Delta(f^n (z))} \le \sum_{z \in \cl{F}(2)} \frac{(\text{diam} \, U_n(z))^t}{\Delta(f(z)) \ldots \Delta(f^n(z))}.
\end{equation}

Continuing the above procedure assume we  have constructed a family $\cl{F}(k) \subset E_n$ and a subfamily $\cl{G}(k)$, so that the sets $(U_n(z))_{z \in \cl{G}(k)}$ 5-cover $W$ and $$\mathop{\sum}\limits_{z \in \cl{G}(k)} \frac{(\text{diam} \, U_n(z))^t}{\Delta(f^{k+1} (z)) \ldots \Delta(f^n (z))} \le  \mathop{\sum}\limits_{z \in \cl{F}(k)} \frac{(\text{diam} \, U_n(z))^t}{\Delta(f(z)) \ldots \Delta(f^n(z))}.$$
For each $w \in W$, we then take the $k$-th iterate of $w(n-1)$ and denote it by $w(n-k-1)$; this is an $(n-k-1)$-preimage of $w$ in $\Lambda$. Now, if $w(n-k-1) \in V_1$ then it has $d_1$ 1-preimages in $\Lambda$ and to each of these we can apply the same procedure from step $k$.  In this way we obtain by succesive eliminations $d_1$ mutually disjoint families $\cl{F}(k+1, d_1; i), 1 \le i \le d_1$ and  inside each of these a subfamily $\cl{G}(k+1, d_1; i)$ such that
$$
\mathop{\sum}\limits_{z \in \cl{G}(k+1, d_1; i)} \frac{(\text{diam} \, U_n(z))^t}{\Delta(f^{k+1} (z)) \ldots \Delta(f^n (z))} \le  \mathop{\sum}\limits_{z \in \cl{F}(k+1, d_1; i)} \frac{(\text{diam} \, U_n(z))^t}{\Delta(f(z)) \ldots \Delta(f^n(z))}.
$$
The succesive elimination procedure works, since we cannot have two differerent 1-preimages $y$ and $y'$ of $w(n-k-1)$ having $(n-k)$-preimages  $\xi\in \Lambda$ and $ \xi' \in \Lambda$ respectively, such that $\xi \in B_n(z, \vp), \xi' \in B_n(z', \vp)$, for some  $z, z' \in \cl{F}(k+1, d_1; i)$. Indeed,  it would then follow that the family  $\{U_n(z): z\in \cl{F}(k+1, d_1; i)\}$ does not consist of pairwise disjoint sets, which clearly  is a contradiction. Moreover, 
since $V_1 \cap V_2 = \emptyset$, we must have $\cl{F}(k+1, d_1; i) \cap \cl{F}(k+1, d_2; j) = \emptyset$. Hence, there is no repetition of elements, when we consider the union
$$
\cl{F}(k+1):= \mathop{\cup}\limits_{1 \le j \le d_1} \cl{F}(k+1, d_1; j) \cup \mathop{\cup} \limits_{ 1\le j \le d_2} \cl{F}(k+1, d_2; j).$$

Now among the collections $\cl{G}(k+1, d_1; i)$, for $1 \le i \le d_1$, let us consider the one which gives rise to the least sum $\sum_{z \in \cl{G}(k+1, d_1; i)} \frac{(\text{diam} \, U_n(z))^t}{\Delta(f^{k+2}(z)) \ldots \Delta(f^n (z))}$. Denote this minimizing collection by $\cl{G}(k+1, d_1)$. Similarly, we obtain the collection $\cl{G}(k+1, d_2)$. We now have that
\begin{equation}\label{stepk}
\aligned
\mathop{\sum}\limits_{z \in \mathcal{G}(k+1, d_1)} \frac{(\text{diam} \, U_n(z))^t}{\Delta(f^{k+2}(z)) \ldots \Delta(f^n(z))} &\le \mathop{\sum}\limits_{z \in \mathop{\cup}\limits_{1 \le i \le d_1} \cl{G}(k+1, d_1; i)} \frac{(\text{diam} \, U_n(z))^t}{\Delta(f^{k+1} (z)) \ldots \Delta(f^n (z))} \\
&\le \mathop{\sum}\limits_{z \in \mathop{\cup}\limits_{1 \le i \le d_1} \cl{F}(k+1, d_1; i)} \frac{(\text{diam} \, U_n(z))^t}{\Delta(f(z)) \ldots \Delta(f^n(z))}.
\endaligned
\end{equation}
Of course, we can proceed similarly for $\cl{G}(k+1, d_2)$. With $\cl{G}(k+1):= \cl{G}(k+1, d_1) \cup \cl{G}(k+1, d_2)$., it follows  from above that
$$
\mathop{\sum}\limits_{z \in \mathcal{G}(k+1)} \frac{(\text{diam} \, U_n(z))^t}{\Delta(f^{k+2}(z)) \ldots \Delta(f^n(z))} \le \mathop{\sum}\limits_{z \in \mathop{\cup}\limits_{1 \le i \le d_1} \cl{F}(k+1)} \frac{(\text{diam} \, U_n(z))^t}{\Delta(f(z)) \ldots \Delta(f^n(z))}.
$$

Therefore, we obtain by finite induction a union $\cl{F}(n)$ of families in $ E_n$, as well as one particular family $\cl{G}(n)$ such that $\{5U_n(z): z \in \cl{G}(n)\}$ covers the set $W$ and has the property that
$$
\mathop{\sum}\limits_{z \in \cl{G}(n)} (\text{diam} \, U_n(z))^t \le \mathop{\sum}\limits_{z \in \cl{F}(n)} \frac{(\text{diam} \, U_n(z))^t}{\Delta(f(z)) \ldots \Delta(f^n (z))}.
$$
By combining this with the observation in (\ref{P}) at the start of the proof, this shows that \[
\mathop{\sum}\limits_{z \in \cl{G}(n)} (\text{diam} \, U_n(z))^t < 1.\]
Since $\{5U_n(z): z \in \cl{G}(n)\}$ is a covering of the set $W=W^s_r(x) \cap \Lambda$, we can now conclude that 
$$
\delta^{s} (x)  \le t< t_\omega.
$$

In  the more general case in which $\omega$ is  a continuous function on $\Lambda$ with the property that
 $\omega(x) \le \Delta(x)$, for all $ x \in \Lambda$, we proceed as follows.
First note that, by continuity of $\omega$, we have that there exists an increasing, positive function $\rho$ on $\Lambda$ such that  $\rho(\vp)$ decreases to zero for $\vp$ tending  to zero from above, and such that if
$
d(y, z) \le \vp$, then $$|\omega(y) - \omega(z)| \le \rho(\vp)$$
Since if $y \in B_n(z, \vp)$ then $f^i (y) \in B(f^i z, \vp)$, the latter implies that if $y \in B_n(z, \vp)$ then  $|\omega(f^i (y)) - \omega(f^i (z))| \le \rho(\vp)$. Hence, since $\Delta(x) \ge \omega(x)$ for all $ x \in \Lambda$, it follows that for each $ 0 \le i \le n-1$ we have  $$\Delta(f^i (y)) \ge \omega(f^i (y)) \ge \omega(f^i (z)) -\rho(\vp).$$ 
Now in order to proceed, let us define the $\vp$-pressure function $P_\vp$, for some arbitrary potential function $\psi$, by
  \[P_\vp(\psi):= \mathop{\liminf}\limits_{n \to \infty} \frac 1n \log \inf \left\{\mathop{\sum}\limits_{x\in E}\exp\left(\sum_{k=0}^{n-1}\psi(f^{k}(x))\right): E \mbox{ is a } (n, \vp)\mbox{-spanning set for} \ \Lambda \right\}\]
and let $t_\vp$ denote the unique zero of $P_\vp(t\Phi^s - \log (\omega- \rho(\vp)))$.
Then let $t > t_\vp$ be fixed and note that the above proof goes through in the same way if in the sums  appearing there, we replace the function $\Delta$ by the function $\omega - \rho(\vp)$.  
Indeed, this follows since for all $0 \le i \le n-1$ we have that  $\Delta(f^iy) \ge \omega(f^i(y)) \ge \omega(f^i(z)) -\rho(\vp)$, for each $y \in B_n(z, \vp)$ and
for some arbitrary fixed element $z$ contained in some minimal $(n, \vp)$-spanning set $E_n$ for $\Lambda$.
In this way,  the above inductive procedure gives rise to a family $\cl{F}(n) \subset E_n$ and to a particular family $\cl{G}(n)$ such that $\{5U_n(z): z \in \cl{G}(n)\}$ covers the set $W$ and such that
$$
\mathop{\sum}\limits_{z \in \cl{G}(n)} (\text{diam} \, U_n(z))^t \le \mathop{\sum}\limits_{z \in \cl{F}(n)} \frac{(\text{diam} \, U_n(z))^t}{(\omega(f(z)) - \rho(\vp)) \ldots (\omega(f^n(z)) - \rho(\vp))} < 1.
$$
Now, for $\eta >0$ sufficiently small and $0 < \vp < \eta$,  let  $\tau_{\vp, \eta}$ refer to the unique zero of the pressure
 function $P_\vp(t\Phi^s - \log (\omega - \rho(\eta)))$ and let $\tau_\eta$ denote the unique zero of the pressure function $P(t\Phi^s - \log (\omega - \rho(\eta)))$.  Since $\mathop{\lim}\limits_{\vp \to 0} P_\vp(\psi) = P(\psi)$ for each continuous function $\psi$, it follows that $\lim_{\vp \to 0}\tau_{\vp, \eta} = \tau_\eta$. 
On the other hand, note that for $0 < \vp < \eta$ we have that $\rho(\vp) < \rho(\eta)$ and therefore, $t\Phi^s - \log(\omega - \rho(\vp)) \le t\Phi^s - \log(\omega - \rho(\eta))$. This implies that  $\tau_\vp \le \tau_{\vp, \eta}$. Now, consider some arbitrary fixed $t > \tau_\eta$. For $\vp>0$ sufficiently small, we then have that $t > \tau_{\vp, \eta} \ge \tau_\vp$. Hence, from the above we have that for  $t$ in this range and for $n$ sufficiently large,  there exists a cover $\{5U_n(z):  z \in \cl{G}(n)\}$ of $W$ such that $$\mathop{\sum}\limits_{z \in \cl{G}(n)} (\text{diam} \, U_n(z))^t < 1.$$ This shows that $t \ge \dim_{H}(W)$ and therefore, since $t>\tau_\eta$ was chosen to be arbitrary,  it follows that $\tau_\eta \ge \dim_{H}(W)$. 
Finally, observe that the continuity of the pressure function implies that $\lim_{\eta \to 0 }\tau_\eta = t_\omega$, and this then allows to deduce the desired inequality $$\dim_{H}(W) \le t_\omega.$$

\, \, $\mbox{ }$ \hfill $\Box$

\section{Proofs of Proposition \ref{max} and Proposition \ref{upperboxdim}.} 

%
%
%
\noindent {\em Proof of Proposition \ref{max}.}

%
Recall that here we assume that $d$ is the minimal value of $\Delta$ on $\Lambda$. Then note that, since $\Delta$ is upper semi-continuous on $\Lambda$ and takes on only integer values, it follows that if $\Delta(x) = d$ for some $x \in \Lambda$, then we have  that the preimage counting function $\Delta$ must be equal to $d$ on some open neighbourhood of $x$. This implies that
the set $$A:= \{x \in \Lambda: \Delta(x) = d\}$$ has to  be open in $\Lambda$. In order to show  that $A$ is dense in $\Lambda$, assume to the contrary that there exists  a non-empty open set $V \subset \Lambda$ such that $\Delta(x) \ge d+1$, for all $x \in V$. In this situation we can then find a Lipschitz continuous function $\psi$ on $\Lambda$ such that $d \le \psi(x) \le \Delta(x)$, for all $x \in \Lambda$, and such that  $\psi \equiv d+1$ on some open ball contained in $V$. 

Now note that Theorem \ref{upper} implies  that $t_\psi \ge \delta^s(x)$, for all $x \in \Lambda$. Also, since  $ \psi(x)\geq d$ for all $x \in \Lambda$, we have that $t_\psi \le t_d$. Therefore,  if for some $x \in \Lambda$  we have that $t_d = \delta^s(x)$, then it follows  that $$t_d = t_\psi = \delta^s(x).$$
Let us now consider the unique equilibrium measure $\mu_\psi$ for the H\"older continuous potential $t_d\Phi^s - \log \psi$ (note that existence and uniqueness of $\mu_\psi$ is guaranteed, since $f$ is hyperbolic on $\Lambda$ (see \cite{KH} and \cite{M-DCDS06})).  Also, since $\mu_\psi$ is a $f$-invariant probability measure for which the Variational Principle holds for the potential $t_d \Phi^s - \log d$, we have that
 $$0 = P(t_d\Phi^s - \log d) = P(t_d \Phi^s - \log \psi) = h_{\mu_\psi} + \int_\Lambda (t_d \Phi^s - \log \psi) \ \d\mu_\psi \ge h_{\mu_\psi} + \int_\Lambda (t_d \Phi^s - \log d) \ \d\mu_\psi.$$ 
This shows that  $$\int_\Lambda \log \psi  \ \d\mu_\psi \le \int_\Lambda \log d \ \d\mu_\psi.$$
However, recall that  $\log \psi(y) > \log d$, for all $y$ in some open ball contained in $V$. Moreover,  since $\mu_\psi$ is an equilibrium measure, we have that   $\mu_\psi$  is positive on Bowen balls  and hence,
it is positive  on any open set in $\Lambda$. Clearly, this gives a contradiction and therefore, it follows that $\Delta \equiv d$ on a dense open set in $\Lambda$. 

In order to show that if $\Delta \equiv d$ on an open dense set then it follows that $\delta^s(y) = t_d$ for all $y \in \Lambda$, we define the set
\[ A_n:= \{x \in \Lambda: x \mbox{ has precisely $d^n$ $n$-preimages $y_{i}$ and } \Delta(f^j (y_{i})) = d, \mbox{ for all } 0 \le j \le n \mbox{ and $1\leq i \leq d^{n}$}\}.\]
The aim is to show that  $A_n$ is open and dense in $\Lambda$, for each $n\in \N$. For this, we first show that  $A_1$ is open in $\Lambda$.  By definition, we have that if $x \in A_1$ then $x\in A$ and hence, $x$ has precisely $d$ 1-preimages $x_1, \ldots, x_d \in A$. Now, let $y$ be a point  close to $x$. Since $A$ is open, we can assume without loss of generality that $y\in A$ and hence, $y$
has precisely $d$ preimages $y_{1},...,y_{d} \in \Lambda$. Since $d$ is the least value $\Delta$ can attain on $\Lambda$ and since $f$ has no critical points in $\Lambda$, we have that each of the $y_{i}$ is close to one of the $x_{j}$. Since $A$ is open and since the $x_{j}$ are contained in $A$, it follows that  $y_i \in A$, for all $1 \le i \le d$ provided $y$ is close enough to $x$. This shows that $y \in A_{1}$ and hence  it follows that  $A_1$ is open in $\Lambda$. 

In order to show that $A_1$ is dense in $\Lambda$, consider some open set $V$ in $\Lambda$. Since $A$ is dense in $\Lambda$,  there exists some point $y \in A \cap V$, which must have precisely  $d$ 1-preimages $y_1, \ldots, y_d \in \Lambda$. 
Now, let $B \subset A$ be a small ball centred at $y$. For each $1 \le i \le d$, choose a sufficiently small  ball $B_i$ centred at $y_i$ such that the resulting family of balls is pairwise disjoint and  such that $f$ is injective on $B_i$ and on $B \subset f(B_i)$.  The aim is to show that $B \cap f(B_i \cap A)$ is open and dense in $B$. Indeed, if  $z \in B \cap f(B_i \cap A)$, then $z$ has a 1-preimage $z_i\in B_i \cap A$. 
Now, if $z'$ is close enough to $z$, then $z'$ belongs to $A$ and hence, $z'$ has a 1-preimage $z'_i \in B_i$ which lies close to $z_i$.  Since $z_i \in A$ and since $A$ is open,  it follows that $z'_i \in A$. This gives that  $B \cap f(B_i\cap A)$ is open in $B$. Also, if there were a non-empty open set $ B' \subset B$ such that $B' \cap f(B_i \cap A) = \emptyset$, then $B_i \cap f^{-1}(B') $ would be  open and non-empty. Clearly, this contradicts the fact that $A$ dense in $\Lambda$. This shows that $B \cap f(B_i\cap A) $ must be open and dense in $B$, for all $1 \le i \le d$. Since a finite intersection of open and dense subsets is again open and dense, it now follows that $A_1$ has to be open and dense in $\Lambda$. 

Clearly, the same methods as in the previous argument  can  be used to prove by way of induction that $A_n$ is open and dense
 in $\Lambda$, for each $n\in \N$.
Therefore, we now have that for each $n\in \N$, there exists an open dense set $A_n$ such that for every $y \in A_n$ there exist exactly $d^n$ $n$-preimages $y_{1},...,y_{d^{n}} \in \Lambda$ of $y$ such that $\Delta(f^i y_{j}) = d$, for each $0 \le i \le n$ and $1\leq j \leq d^{n}$ . \ 
This shows that in the proofs of Theorem \ref{upper} and the main theorem of \cite{MU-BLMS} one can work exclusively  with points from $\bigcup_{n\in \N}A_n$. Indeed since $A_n$ is open and dense in $\Lambda$, it follows that for every $z \in E_n$ we can take a point $z'$ sufficiently close to $z$ such that $f^n(z') \in A_n$; thus we obtain a set $E_n'$ with the same number of elements as $E_n$ which is again $(n, \vp)$-spanning and can be used in the condition on the pressure in order to obtain good covers of $W^s_r(x) \cap \Lambda$. Then all the iterates up to order $n$ of any $z' \in E_n'$ will have exactly $d$ 1-preimages in $\Lambda$ and then we obtain  $\delta^s(x)=t_{d}$, for all $x \in \Lambda$.

\, \, $\mbox{ }$ \hfill $\Box$

\

\noindent {\em Proof of Proposition \ref{upperboxdim}.}

%
%
%
(a) In the sequel let $x \in \Lambda$ be fixed and put $W:= W^s_r(x) \cap \Lambda$. As in the proof of Theorem \ref{upper}, for each $\vp>0$ sufficiently small,  there exists $n_0\in \N$ and a minimal $(n_0, \vp)$-spanning set $E_{n_0}$ for $\Lambda$ such that for each $t>t_\omega$ sufficiently large we have, for some fixed $\beta>0$,  
\begin{equation}\label{delta}
\sum_{z \in E_{n_0}} \frac{|Df_s^{n_0}(z)|^t}{\omega(f(z) \ldots \omega(f^{n_0} z)} < e^{-\beta n_0} < 1/2.
\end{equation}

Let us assume $E_{n_0} =: \{e_1, \ldots, e_{m_0}\}$. As before,  define $U_n(z):=f^n(B_n(z, \vp)) \cap W^s_r(x)$,  for $n\in \N$ and $ z \in \Lambda$.
The aim is to construct a covering of $W$ which consists of sets of comparable diameter. For this, let  $\{|Df_s^{n_0}(z)|: z \in E_{n_0}\}=:\{\delta_1, \ldots, \delta_{m_{0}}\}$ and then define 
for  $n \in \N$ the value $\chi(n)$ by  
$$\chi(n):= \inf \left\{ \prod_{i=1}^{n}\delta_{j_i} : 1 \le j_i \le m_{0} \right\}.$$
Now,  for each $w \in \Lambda$ and for each $n n_0$-preimage $w(-n n_0) \in \Lambda$ of $w$, we have that  $f^{jn_0}(w(-nn_0)) \in B_{n_0}(z_j, \vp)$, for each $0 \le j \le n-1$. From this we deduce that  $|Df_s^{nn_0}(w(-nn_0))| \ge \chi(n)$. Next 
observe that in general, given any full prehistory $\hat w=(w, w_{-1}, \ldots) \in \hat \Lambda$ of some element $w\in \Lambda$,  there exists $k(\hat w, n) \in \N$ such  that $|Df_s^{k(\hat w, n) n_0}(w_{-k(\hat w, n) n_0})| $ is comparable to $\chi(n)$, that is, 
$$ C_{0}^{-1} \cdot \chi(n) < |Df_s^{k(\hat w, n) n_0}(w_{-k(\hat w, n) n_0})| < C_{0} \cdot \chi(n),$$
where we have put $C_{0}:= \sup_{z \in \Lambda} \cdot |Df_s^{n_0}(z)| $.
This shows that for $w \in W$ we have that the  diameter $\text{diam} \, U_{ k(\hat w, n) n_0}(w_{-k(\hat w, n) n_0 })$ is comparable to $\chi(n)$, where the comparability constant does  depend neither on $w$ nor on $n$. Hence,  the sets $U_{ k(\hat w, n)n_0}(w_{- k(\hat w, n) n_0})$ provide a covering of $W$ and their diameters are all of size comparable to $\chi(n)$.
For later use, let us remark that one can choose  a point $z_{k(\hat w, n)}(\hat w) \in E_{n_0}$ such that $w_{-n_0 k(\hat w, n)} \in B_{n_0}(z_{k(\hat w, n)}(\hat w), \vp)$ and similarly, points $z_{k(\hat w, n) - j}(\hat w) \in E_{n_0}$ such that $f^{n_0j}(w_{-n_0 k(\hat w, n)}) \in B_{n_0}(z_{k(\hat w, n) -j}(\hat w), \vp)$, for each $ 1 \le j < k(\hat w, n)$. \
 Then recalling that $E_{n_0} =: \{e_1, \ldots, e_{m_0}\}$, the inequality in (\ref{delta}) reads:

$$
\sum_{i=1}^{m_{0}} \frac{\delta_i^t}{\omega(f (e_i)) \ldots \omega(f^{n_0}  (e_i))} < \frac 12.
$$
By raising both sides of this inequality to the power $p\in \N$ and then summing  over $p$, we obtain
\begin{equation}\label{u}
\aligned
\sum_{p \in \N} & \left(\sum_{i=1}^{m_{0}} \frac{\delta_i^t}{\omega(f(e_i)) \ldots \omega(f^{n_0}  (e_i))}\right)^p = \\
&=\sum_{p\in \N} \sum_{(i_{1},...,i_{p}) \in \{1,...,m_{0}\}^{p}} \frac{\delta_{i_1}^t \ldots \delta_{i_p}^t}{\left(\omega(f (e_{i_1})) \ldots \omega(f^{n_0} (e_{i_1}))\right)\cdot \ldots \cdot \left(\omega(f (e_{i_p})) \ldots \omega(f^{n_0} (e_{i_p}))\right)} < 1.
\endaligned
\end{equation}

Let us now again consider  some point $w \in \Lambda$ and its full prehistory $\hat w=(w, w_{-1}, \ldots) \in \hat \Lambda$. By the  above, we then have that the orbit of $w_{- k(\hat w, n) n_0}$ under the map $f^{ k(\hat w, n) n_0}$ is shadowed by the consecutive linking of the $n_0$-orbits of $k(\hat w, n)$ points from $E_{n_0}$. Then, the summand of the corresponding sum,  associated with this orbit, is  of the form
\begin{equation}\label{t}
\frac{  \left(\text{diam} \, U_{ k(\hat w, n) n_0}(w_{- k(\hat w, n)n_0})\right)^{t}}{\left(\omega(z_{k(\hat w, n)}(\hat w)) \ldots \omega(f^{n_0}(z_{k(\hat w, n)}(\hat w)))\right) \ldots \left(\omega(z_{1}(\hat w)) \ldots \omega(f^{n_0}(z_1(\hat w)))\right)}.
\end{equation}
We can now use the procedure of successive partial minimisation and elimination, which we used in the proof of Theorem 
\ref{upper}, and this then leads to a covering of $W^s_r(x) \cap \Lambda$ consisting of sets of diameter comparable to $\chi(n)$. Indeed,  as in the proof of Theorem \ref{upper}, here we use the fact that the denominators of the terms in (\ref{u}) are products of evaluations of $\omega$ along the forward orbit of the preimages.  
In this way we obtain a sum with summands of the form as in (\ref{t}), which is smaller or equal than the sum in (\ref{u}). To this sum we can apply  the repeated partial minimisation procedure as in the proof of Theorem 
\ref{upper}, in order to extract a subcover $\mathcal{V}$ such that in the associated sum the denominators are successively eliminated, that is,  we arrive at the inequality
\[ \sum_{U \in  \mathcal{V}}\left(\text{diam} \, U\right)^{t}<1.\]

From this it clearly follows that  $$\beta^s (y) \le t_\omega, \mbox{ for all } y \in W^s_r(x) \cap \Lambda.$$
 
\

 (b) \ The aim is to show that the stable upper box-counting dimension $\beta^s$ is constant on $\Lambda$. 
For this note that, since $f$ is transitive on $\Lambda$,  there exists a point $x\in \Lambda$ whose set of preimages is dense in $\Lambda$.  Therefore,  if $y\in \Lambda$ is some fixed arbitrary point  and if $\vp>0$, then there exists  some $n$-preimage $x_{-n}$ of $x$ such that $d(x_{-n}, y) = \vp$, for some $n\in \N$. 

Then notice that  the local product structure (see \cite{KH}) implies that if for some $z \in \Lambda$ the  local unstable manifold $W^u_r(\hat z)$ intersects $W^s_r(y)$, then it will intersect  $W^s_r(x_{-n})$ at a unique point contained in $\Lambda$. Likewise, any local unstable manifold which intersects $W^s_r(x_{-n})$ will also  intersect $W^s_r(y)$ in a point from $\Lambda$. Note that if $W^s_r(y)\cap \Lambda$ is covered by balls $U \in \mathcal{U}$ of radius $\vp>0$, then the set $W^s_r(x_{-n}) \cap \Lambda$ is covered by the same number of balls of radius at most $C' \vp$, for some fixed constant $C'>0$. This follows, since the intersection  $W^s_r(x_{-n}) \cap \bigcup_{\hat z \in \hat \Lambda, z \in U}W^u_r(\hat z)$
 is contained in a ball of radius $C' \vp$, which follows since $d(x_{-n}, y) = \vp$ and since the inclination of local unstable manifolds with respect to $W^s_r(y)$ is bounded from below, a consequence of the uniform hyperbolicity of $f$ on $\Lambda$. 

Also, if we cover $W^s_r(x_{-n}) \cap \Lambda$ with balls  of radius $\vp$, then we can consider all local unstable manifolds through the points of each of these balls to obtain balls  of radius at most $C' \vp$  which are contained in these balls in $W^s_r(y)$.  However, by setting $\vp' := \vp |Df_s(x_{-n})|^n$ for $\vp>0$ sufficiently small, we have that every covering by balls of radius $\vp'$ of $W^s_r(x) \cap \Lambda$ gives a covering by balls of radius $ \vp$ of $W^s_r(x) \cap \Lambda$. Therefore $$\beta^s(y) = \beta^s(x), \mbox{ for all } y \in \Lambda$$ 
and therefore, it follows that the stable upper box dimension is constant throughout $\Lambda$.   
\, \, $\mbox{ }$ \hfill $\Box$

{\em Remark.} Let us assume for a moment that  the following condition is satisfied: if $D$ is the maximum possible value of $\Delta$ on $\Lambda$, then for each $ 1 \le i \le D-1$ the sets $\Lambda_{i}:= \{x \in \Lambda: \Delta(x) \le i\}$ have their respective closure contained in $\Lambda_{i+1}$. 
Note that, by the upper semi-continuity of $\Delta$ on $\Lambda$, we have that the set $\Lambda_D:=\{x \in \Lambda: \Delta(x) =D\}$ is closed in $\Lambda$. Also, the upper semi-continuity of $\Delta$ implies that  $\Lambda_i$  is  open  in $\Lambda$, for each $1 \le i \le D-2$. Due to our assumption here, it is possible to fix some neighbourhood $\Lambda_i(\vp)$ of $\bar \Lambda_i$ such that $\Lambda_i(\vp) \subset \Lambda_{i+1}$, for each $1 \le i \le D-2$. Also, let us fix some neighbourhood $\Lambda_{D-1}(\vp)$ of the closure of $\Lambda_{D-1}$.
Then define $K_0:= \Lambda_D\setminus \Lambda_{D-1}(\vp)$, $K_1:= \bar \Lambda_{D-1} \setminus \Lambda_{D-2}(\vp), K_2:= \bar \Lambda_{D-2} \setminus \Lambda_{D-3}(\vp), \ldots, K_{D-1}:= \bar \Lambda_1$ and 
note that the family $\{K_j: 0 \le j < D\}$ consists of pairwise disjoint compact sets. Therefore, there exists a continuous function $\psi$ on $\Lambda$ such that $\psi(x) = D$ for all $ x \in K_0$, $D-1 \le \psi(x) \le D$ for $x \in \Lambda_{D-1}(\vp) \setminus \bar \Lambda_{D-1}$, $\psi(x) = D-1$ for $x \in K_1$, and $D-2 \le \psi(x) \le D-1$ for $x \in \Lambda_{D-2}(\vp) \setminus \bar \Lambda_{D-2}$, which can be continued until we reach $\Lambda_1$. By construction, we then have  that $\Delta(x) \ge \psi(x)$, for all $ x \in \Lambda$.
By applying Theorem \ref{upper},  it follows that $\delta^s(x) \le
  t_{\psi_\vp}$, for all $ x \in \Lambda$ and $\vp>0$. Also, by choosing $\vp \ge \vp'$ appropriately,  we can assume that $\Lambda_i(\vp') \subset \Lambda_i(\vp)$. Therefore, we have for each  $x \in \Lambda$ that $\psi_\vp(x)$ is increasing, for  $\vp$ tending to  zero. This implies  that there exists  $t_*$ such that $t_{\psi_\vp}$ tends to $ t_*$, for $\vp$ tending to zero,  and therefore, we have that $\delta^s(x) \le t_*$, for each $ x \in \Lambda$.

\

  \section{Two examples}
  
\ \  \noindent  \textbf{Example 1.}  \
We assume that the reader is familiar with the type of horseshoes introduced by Simon and Solomyak in \cite{SS}. They considered  horseshoes with overlaps in $\mathbb R^3$ which are given by a $\mathcal{C}^{1+\epsilon}$-transformation $f$, defined by
$$
f(x, y, z) := (\gamma(x, z), \eta(y, z), \psi(z)), \mbox{ for all }  (x, y, z) \in [0,1] \times [0,1] \times \mathcal{I},
$$
 where $\mathcal{I} := \bigcup_{i=1}^{m}\in I_i$ denotes the  union of $m$ compact pairwise disjoint intervals $I_1, \ldots, I_m \subset (0,1)$; we also assume that $m \ge 3$, that $\lambda_1 < |\gamma'_x|, |\eta'_y| < \lambda_2$ for some $0 < \lambda_1 < \lambda_2 < 1/2$, that $|\psi'| > 1$ on $\mathcal{I}$, and that $\psi(I_i) = [0, 1]$, for all $i = 1, \ldots, m$.
The basic set $\Lambda$ of $f$  is defined as before, that is, $\Lambda:= \cap_{n \in \mathbb Z} f^n([0, 1]^3)$. 
Let us now consider  the following perturbations $f_\tau$ of $f$:
\begin{equation}\label{S}
f_\tau(x, y, z) := (\gamma(x, z) + \tau_{i, 1}, \eta(y, z) + \tau_{i, 2}, \psi(z)), \mbox{ for all }  (x, y, z) \in [0, 1] \times [0,1] \times I_i, 1 \le i \le m.
\end{equation}
  We will say that $\tau := (\tau_{1, 1}, \tau_{1, 2}, \ldots, \tau_{m, 1}, \tau_{m, 2})$ is $f$-\textit{admissible} if $f_\tau(\bigcup_{1 \le i \le m} [0, 1]^2 \times I_i) \subset (0, 1)^{2} \times [0, 1]$. It can be checked that  the set of $f$-admissible parameters $\tau$ is a non-empty open subset of $\mathbb R^{2m}$. Also, due to the expansion in the $z$-direction as well as the contractions with respect to the 
 $(x, y)$-coordinates, one can show that $f_\tau$ is hyperbolic on the basic set  $\Lambda_\tau$ associated with $f_\tau$.

As in \cite{SS}, one then verifies that for Lebesgue almost every $f$-admissible $\tau$ we have that the stable dimension of  $\Lambda_\tau$ is given by the maximum of the zeros $s_1$, and $ s_2$ respectively, of certain pressure functions of $\log|\gamma'_x|$, and $\log |\eta'_y|$ respectively, on the symbolic space $\Sigma_{m}$. Let us now assume  that on  $[0, 1] \times \mathcal{I}$ we  have $$|\gamma'_x| = |\eta'_y| \equiv  1/m$$  Then from the proof of Theorem 1 $i)$ of \cite{SS} and the fact that
in this case both zeros $s_1$ and $ s_2$ are  equal to 1,  it follows that the stable dimension of $\Lambda_\tau$ is equal to $1$, for Lebesgue-almost every $f$-admissible $\tau$. 

However, in the above case we have that  the zero $t_{1, \tau}$ of the pressure function
 $t \mapsto P_{f_\tau|_{\Lambda_\tau}}(t\Phi^s_\tau)$ for the stable potential function $\Phi^s_\tau$,  is also equal to $1$. This follows, since  $\Phi^s_\tau \equiv -\log m$ and since  the entropy of $f_\tau|_{\Lambda_\tau}$ is equal to $\log m$, where the latter is due to the fact that the spanning sets of $f|_{\Lambda_\tau}$ are determined only by the dynamics of $\psi$ in the $z$ coordinate; but  this dynamics in the $z$-direction is conjugated to the shift $\sigma_m$ on $\Sigma_m$, since $\psi$ expands $I_i$ onto the whole interval  $[0, 1]$ for each  $i = 1, \ldots, m$. This shows that  $t_{1, \tau} = 1$. Also note  that if $|\gamma'_x| = |\eta'_y| \equiv  1/m$ on  $[0, 1] \times \mathcal{I}$, then  $f_\tau$ is c-conformal.  Therefore, since $\delta^s = t_{1, \tau}$, we can now apply  Corollary \ref{inj}, which then gives that almost every horseshoe $f_\tau$ has an open dense set of points in its associated basic set $\Lambda_\tau$ such that each of these points has
  precisely one $f_\tau$-preimage in $\Lambda_\tau$.

In conclusion, for the above choice of $\gamma, \eta$, we have now shown that Lebesgue-almost every translation $f_\tau$ is close to being a homeomorphism on its associated basic set $\Lambda_\tau$. We summarise this result in the following:

\begin{cor}\label{horseover}
 Let $(f_\tau)_\tau$ denote the family of horsehoes with overlaps given in (\ref{S}), and assume that on  $[0, 1] \times \mathcal{I}$ we have $|\gamma'_x| = |\eta'_y| \equiv  1/m$. Also, let  $\Lambda_\tau:= \cap_{n \in \mathbb Z} f_\tau^n([0, 1]^3)$ denote  the associated basic set of $f_\tau$. Then, for Lebesgue-almost every $f$-admissible parameter $\tau$ there exists an open dense set $A_\tau$ in $\Lambda_\tau$, such that every $x \in A_\tau$ has precisely one $f_\tau$-preimage in $\Lambda_\tau$.
\end{cor}  

\

 \ \ \noindent \textbf{Example 2.} \
In \cite{M-MZ} the first author gave an example of a family of non-linear hyperbolic skew products
 for which the preimage counting function is not constant on their associated basic sets.
Let us first briefly recall the construction of this family. 
For $\alpha \in (0, 1)$, let  $I_1^\alpha, I_2^\alpha \subset I := [0 , 1]$ be two intervals  such that
$I_1^\alpha \subset [\frac 12 - \epsilon(\alpha), \frac
12 + \epsilon(\alpha)]$ and $I_2^\alpha \subset [1-\alpha - \epsilon(\alpha), 1-\alpha + \epsilon(\alpha)]$, for 
some  $0<\epsilon(\alpha) < \alpha^2$ sufficiently small. Let  $g:I_1^\alpha \cup I_2^\alpha \to I$ be a strictly
increasing smooth function with the property
that $g(I_1^\alpha)= g(I_2^\alpha) = I$. Also, assume that there exists a large number $\beta>0$ such that $\beta^2
> g'(x)
> \beta$, for each $ x \in I_1^\alpha \cup I_2^\alpha$. Then there exist intervals
$I_{11}^\alpha, I_{12}^\alpha \subset I_1^\alpha$ and $ I_{21}^\alpha,
I_{22}^\alpha \subset I_2^\alpha$ such that $g(I_{11}^\alpha) =
g(I_{21}^\alpha) = I_1^\alpha$ and $g(I_{12}^\alpha) =
g(I_{22}^\alpha) = I_2^\alpha$. For $J^\alpha:=
I_{11}^\alpha \cup I_{12}^\alpha \cup I_{21}^\alpha \cup
I_{22}^\alpha$ and $J_*^\alpha:= \{x\in J^\alpha: g^i(x) \in
J^\alpha \mbox{ for all } i \ge 0\}$, we then let 
 $f_\alpha: J_*^\alpha \times I
\to J_*^\alpha \times I$ be defined by
\begin{equation}\label{skew-general}
f_\alpha(x, y) := (g(x), h_\alpha(x, y)), \mbox{ where }
h_\alpha(x, y):= \left\{ \begin{array}{ll}
                     \psi_{1, \alpha}(x) + s_{1, \alpha} y, \ x \in I_{11}^\alpha \\
                     \psi_{2, \alpha}(x) + s_{2, \alpha} y, \ x \in I_{21}^\alpha \\

                     \psi_{3, \alpha}(x) - s_{3, \alpha} y, \ x \in I_{12}^\alpha \\
                      s_{4, \alpha} y, \ x \in I_{22}^\alpha,
                     \end{array}
           \right.
\end{equation}
where  $s_{1, \alpha}, ..., s_{4, \alpha} \in (1/2 -\vp_0, 1/2 +\vp_0)$ denote some arbitrary fixed numbers close to $1/2$  and $\psi_{1, \alpha}, \psi_{2,
\alpha}, \psi_{3, \alpha}:I \to \R$ are 
$\mathcal{C}^2$-functions which are $\vp_0$-close (with respect to the
$\mathcal{C}^1$-metric) to the linear functions given by $x \mapsto x$, $x \mapsto
1-x$ and $x \mapsto 1$ respectively.
Let us also use the following shorter notation:
$$h_{x, \alpha}(y):= h_\alpha(x, y), \mbox{for any $(x, y)$ for which this is well-defined}.$$ 
By defining $h^n_{z, \alpha}:= h_{f^{n} (z), \alpha} \circ \ldots \circ h_{z, \alpha}$ for each $ n \ge 0$, the basic set $\Lambda_\alpha$ of the above system is  given by \\ 
$$
\Lambda_\alpha = \bigcup_{x \in J_*^\alpha}\bigcap_{n \ge 0} \bigcup_{z \in g^{-n}(x) \cap J_*^\alpha} h_{z, \alpha}^n(I).$$
In \cite{M-MZ} it was shown that for $\alpha$ small enough, the map $f_\alpha$ is a hyperbolic endomorphism on $\Lambda_\alpha$ and  that there exist two infinite point sets $A_\alpha, B_\alpha \subset \Lambda_\alpha$, which are both not dense in $\Lambda_{\alpha}$, such that for  each point in $A_\alpha$ the preimage counting function $\Delta$ is constant equal to $1$, whereas for points in $B_\alpha$ we have that $\Delta$ is constant equal to $2$. Since $\Delta$ is upper semi-continuous, it follows that $A_\alpha$ is open in $\Lambda_\alpha$ and that  $B_\alpha$ is closed in $\Lambda_\alpha$. 
We want now to obtain an estimate on $\delta^s$ by using Theorem \ref{upper}. 
There exists a non-constant continuous function $\omega$ such that  $1 \le \omega(x) \le 2$ for all $ x \in \Lambda_\alpha$, and  $\omega(x)= 1$ for all  $x \in A_\alpha$. Hence we have $\omega(x) \le \Delta(x), x \in \Lambda_\alpha$.  The function $\omega$ may be viewed as an "approximation" of the function $\chi_{B_\alpha} +1$. 
Then an application of Theorem \ref{upper} gives the following estimate:  $$\delta^s(y) \le t_\omega,  \mbox{ for all } \ y \in \Lambda_\alpha$$ 
Similarly as before consider an increasing sequence of non-constant continuous functions $\omega_m$ such that $\omega_m \equiv 1$ on $A_\alpha$,  and $1 \le \omega_m \le 2$ on $\Lambda_\alpha$. For each member of this sequence we can now argue as before. This leads to improvements of the upper bounds for $\delta^s(y)$. Namely,  with  $t_{\omega_m}$ referring to the zero of the pressure function associated with  $\omega_m$, we have that  $t_{\omega_m}$ are decreasing when $m \to \infty$ and $$\delta^s(y) \le t_{\omega_m}, \mbox{ for all } y \in \Lambda_\alpha, m \in \N.$$ In particular,  by applying Proposition \ref{upperboxdim}, we obtain that the stable upper box dimension $\beta^s$ is constant on $\Lambda_\alpha$ and bounded above by  $t_{\omega_m}$, for each $m \in \N$.

$\hfill\square$

\textbf{Acknowledgements:} 

\

\

Eugen Mihailescu, \ \ \ \ \  Eugen.Mihailescu\@@imar.ro

Institute of Mathematics of the Romanian Academy,
P. O. Box 1-764,  RO 014700,  

Bucharest, Romania.

\

Bernd Stratmann, \ \ \ \ \ bos@math.uni-bremen.de

Fachbereich 3 - Mathematik, Universit\"at Bremen, 

Bremen, Germany.


\begin{thebibliography}{99}

\bibitem{Ba}
L. Barreira, \emph{Dimension and recurrence in hyperbolic dynamics}, Progress in Mathematics vol. 272, Birkhauser, 2008.

\bibitem{Bot}
H. G. Bothe,  Shift spaces and attractors in noninvertible horseshoes, Fund. Math. \textbf{152}, 1997, 267-289.

\bibitem{Bo}
R. \ Bowen, \emph{Equilibrium states and the ergodic theory of Anosov
diffeomorphisms}, Lecture Notes in Mathematics, 470, Springer 1975.

\bibitem{Bo-dim}
R. Bowen, Hausdorff dimension of quasicircles, Publ. Math. Inst. Hautes Etudes Sci. \textbf{50}, 1979, 11-25.

\bibitem{E}
G. Edgar, Fractal dimension of self-similar sets: some examples, Suppl. Rend. Circ. Mat. Palermo Ser. II, \textbf{28}, 1992, 341-358.

\bibitem{Fa}
K. Falconer, The Hausdorff dimension of some fractals and attractors of overlapping construction, J. Stat. Phys. \textbf{47}, 1987, 123-132.

\bibitem{KH}
A.\ Katok and B.\ Hasselblatt, \emph{Introduction to the Modern Theory
of Dynamical Systems}, Cambridge Univ. Press, London-New York,
1995.

\bibitem{MM}
A. Manning and H. McCluskey, Hausdorff dimension for horseshoes,  Ergodic Th. and Dynam. Systems
\textbf{3},  1983, 251-260.

\bibitem{Mat}
P. Mattila, \emph{Geometry of sets and measures in euclidean spaces. Fractals and Rectifiability}, Cambridge University Press, 1995.

\bibitem{M-MZ}
E. \ Mihailescu, Unstable directions and fractal dimensions for a family of skew products with overlaps, Math. Zeitschrift, \textbf{269}, 2011, 733-750.


\bibitem{M-DCDS06}
E. \ Mihailescu, Unstable manifolds and Holder structures
associated with noninvertible maps, Discrete and Cont. Dynam.
Syst. \textbf{14}, 3, 2006, 419-446.

\bibitem{MU-BLMS}
E. \ Mihailescu and M. \ Urba\'nski, Relations between stable dimension and the preimage counting
function on basic sets with overlaps, Bull. London Math. Soc. \textbf{42}, 2010, 15-27.






\bibitem{PS} Y. \ Peres and B. \  Solomyak, Problems on self-similar sets and self-affine sets: an update, Progress
in Probability \textbf{46}, 2000,  95-106.

\bibitem{P}
Y. \ Pesin, \emph{Dimension theory in dynamical systems}, Chicago Lectures in Math. Series, 1997.

\bibitem{Pr}
F. Przytycki,  Anosov endomorphisms, Studia Mathematica \textbf{58}, 1976,  249-285.

\bibitem{PU-89}
F. Przytycki and M. Urba\'nski, On Hausdorff dimension of some fractal sets, Studia Math. \textbf{93}, 1989, 155-186.


\bibitem{Ru-carte89}
D. \ Ruelle, \emph{Elements of differentiable dynamics and bifurcation
theory}, Academic Press, New York, 1989.

\bibitem{Ru-carte78}
D. \ Ruelle, \emph{Thermodynamic formalism}, Addison-Wesley, Reading, MA,
1978.

\bibitem{Ru-dim}
D. Ruelle, Repellers for real analytic maps, Ergodic Th. and  Dynam. Systems \textbf{2}, 1982, 99-107.

\bibitem{S}
A. Schief, Separation properties for self-similar sets, Proceed. Amer. Math. Soc., \textbf{122}, 1994, 111-115.

\bibitem{Si}
K. Simon, Hausdorff dimension for non-invertible maps, Ergodic Th. and Dynam. Syst. \textbf{13}, 1993, 199-212.

\bibitem{SS}
K. Simon and B. Solomyak, Hausdorff dimension for horseshoes in $\mathbb R^3$, Ergod. Th. and Dynam. Syst., \textbf{19}, 1999, 1343-1363.

\bibitem{So}
B. Solomyak, Measure and dimension for some fractal families, Math. Proc. Cambridge Phil. Soc. \textbf{124}, 1998, 531-546.


\end{thebibliography}
\end{document}